\documentclass[no preprint line,12pt]{elsarticle}

\usepackage{amssymb}
\usepackage[english]{babel}
\usepackage{amsmath}
\usepackage{array}
\usepackage{caption}
\usepackage{lineno}
\usepackage{float}

\usepackage{graphicx}
\usepackage{color}
\usepackage{vcell}

\usepackage{epstopdf}

\usepackage{subfigure}

\usepackage{algorithm}
\usepackage{algpseudocode}
\usepackage{pifont}
\usepackage{setspace}

\def\be#1\ee{\begin{equation}#1\end{equation}}
\newcommand{\bea}{\begin{eqnarray}}
\newcommand{\eea}{\end{eqnarray}}
\newcommand{\beas}{\begin{eqnarray*}}
\newcommand{\eeas}{\end{eqnarray*}}
\newtheorem{remark}{Remark}

\newcommand{\la}{\langle}
\newcommand{\ra}{\rangle}

\newcommand{\calL}{\mathcal{L}}

\newcommand{\calR}{\mathcal{R}}

\renewcommand{\Re}{\mathfrak{Re} \,}

\def\bI{\mathbf{I}}

\def\bM{\mathbf{M}}

\def\bS{\mathbf{S}}
\def\bT{\mathbf{T}}

\def\bV{\mathbf{V}}

\def\bb{\mathbf{b}}

\def\bee{\mathbf{e}}

\def\bv{\mathbf{v}}
\def\bw{\mathbf{w}}

\def\by{\mathbf{y}}

\def\bsT{\mathsf{T}}

\newcommand{\bphi}{{\boldsymbol{\phi}}}

\newcommand{\CC}{{\mathbb{C}}}

\newcommand{\calQ}{\mathcal{Q}}
\newcommand{\calW}{\mathcal{W}}

\journal{Journal of Computational Physics}

\begin{document}

\begin{frontmatter}



\title{Lippmann-Schwinger-Lanczos algorithm for  inverse scattering problems with unknown reflectivity and loss distributions: One-dimensional Case} 

\author[label1]{J\"orn Zimmerling} 
\author[label2]{Mikhail Zaslavsky}
\author[label3]{Alexander V. Mamonov}
\author[label4]{Vladimir Druskin}
\author[label2]{Anarzhan Abilgazy}
\affiliation[label1]{organization={Uppsala Universitet},
            addressline={Lägerhyddsvägen 2}, 
            city={Uppsala},
            postcode={75237}, 
            state={Uppsala},
            country={Sweden}}
\affiliation[label2]{organization={Southern Methodist University},
            addressline={3100 Dyer St}, 
            city={Dallas},
            postcode={75205}, 
            state={TX},
            country={USA}}
\affiliation[label3]{organization={University of Houston},
            addressline={3551 Cullen Blvd}, 
            city={Houston},
            postcode={77204}, 
            state={TX},
            country={USA}}
\affiliation[label4]{organization={Worcester Polytechnic Institute},
            addressline={100 Institute Road, WPI - Stratton Hall}, 
            city={Worcester},
            postcode={01609}, 
            state={MA},
            country={USA}}

\begin{abstract}
We consider one-dimensional inverse scattering in attenuating media where both the reflectivity and loss distributions are unknown. Mathematically, this corresponds to recovering the coefficients of a damped wave operator, or equivalently, a quadratic operator pencil in the frequency domain. 

The Lippmann-Schwinger equation maps the unknown reflectivity and loss distribution to the measured scattered data. This mapping is nonlinear, as it requires knowledge of the internal wavefield, which itself depends on the reflectivity and loss distribution. The Lippmann-Schwinger-Lanczos method addresses this nonlinearity by approximating the internal solutions through the lifting of states from a reduced-order model constructed directly from the measured data.

In this work, we extend the method to dissipative problems, enabling the approximation of internal partial differential equation (PDE) solutions in media with both reflectivity and loss distributions. We present two complementary constructions of such internal solutions: one based on spectral data and another on frequency-domain measurements over a finite interval. This development establishes a direct link between data-driven reduced-order models for inverse problems and port-Hamiltonian dynamical systems, with reduced models obtained either from the associated spectral measure or via rational approximation. Compared to the Born approximation, which replaces the internal field with the background field, our approach yields more accurate internal reconstructions and enables faster and more robust recovery of the contrast as evidenced by our numerical experiments.

\end{abstract}


\begin{graphicalabstract}
\end{graphicalabstract}

\begin{highlights}
\item Extension of the Lippmann-Schwinger-Lanczos framework to inverse scattering with attenuation, enabling the approximation of internal PDE solutions from boundary measurements.

\item Establishing a link between data-driven reduced-order models and port-Hamiltonian dynamical systems, with reduced models constructed directly from the associated spectral measure and through rational approximation.
\end{highlights}

\begin{keyword}
Inverse Scattering \sep Data-Driven \sep Passive \sep  Lippmann-Schwinger


\MSC 37N30 \sep 65N21 \sep 65L09 \sep 86A22

\end{keyword}

\end{frontmatter}



\section{Introduction}\label{Introduction}

We study one-dimensional inverse scattering problems for wave equations in attenuating (damped) media, where both reflectivity (or wave speed) and loss profiles are unknown. Such problems, for example, arise in underground sensing applications  such as a ground penetrating radar. Mathematically, this corresponds to recovering the coefficients of a second order differential  operator that appears in both the first and second time derivatives, or equivalently, in a quadratic polynomial in the frequency domain \cite{ButerinYurko}. Such problems arise in both electromagnetic and visco-acoustic applications. In the electromagnetic case, they correspond to Maxwell’s equations with both displacement and electro-conductivity currents.

The data-driven reduced-order model (ROM) approach has recently gained prominence as a powerful tool for solving inverse scattering problems when data depend nonlinearly on unknown coefficients \cite{borcea2005continuum, DMTZ16,NormG,TristanSchr,borcea2020reduced,borcea2023data,borcea2023waveform1,borcea2023waveform2,borcea2025reduced,tataris2025inverse}. In such settings, linearization methods such as the Born approximation fail. The generic workflow of the ROM approach can be summarized as follows. First, the measured data (transfer function) are lifted to a ROM, i.e. a dynamical system that reproduces the data and has properties similar to coarse discretization of the underlying PDE. This ROM can then be embedded into the PDE and interpreted as a finite-volume (FV) approximation, yielding estimates of the PDE coefficients from the discrete FV coefficients. In \cite{LossyOneDimensional}, the authors, jointly with Borcea, extended this framework to lossy problems. A key difficulty in that extension was the contradiction between the uniqueness of the discrete transmission-line realization and the inherent non-uniqueness of the continuum problem. This limitation reflects an intrinsic property of dynamical systems: they cannot be transformed into Stieltjes problems with discrete spectra \cite{Zimmerling23}. To overcome this, constrained optimization was introduced in \cite{LossyOneDimensional}, but at the cost of multiple forward solves.

In this work, we propose a frequency-domain variant of the Lippmann-Schwinger-Lanczos (LSL) algorithm \cite{DMZ21,baker25,BorceaZimmerling,Abilgazy}. The method employs a data-driven ROM to approximate the internal PDE solution inside the medium, away from the measurement points. This approximate internal solution is then substituted into the Lippmann-Schwinger integral equation, reducing the problem to a linear one. Such linearization enables a straightforward incorporation of solution constraints. Compared to \cite{Zimmerling23}, our approach requires fewer forward solves for complex lossy media and, importantly, is naturally extensible to higher spatial dimensions.

Unlike the lossless case, inverse scattering in the presence of both propagation and attenuation is intrinsically non-unique. A passive problem with losses cannot be distinguished from a lossless problem in unbounded domains, since both can be represented as Stieltjes functions with nonnegative spectral measures \cite{Zimmerling23}. This ambiguity is more pronounced in the time domain, where finite-time measurements cannot separate losses from long-lived trapped modes, such as those in photonic crystals. For this reason, time-domain variants of the LSL algorithm that use propagator-based ROMs to compute data-driven internal solutions are not applicable in the present setting. By contrast, in the frequency domain the non-uniqueness can be mitigated by assuming that the true problem and the perturbed problem share the same high-frequency asymptotics. This condition can be enforced by constraining the perturbation in the Lippmann-Schwinger equation, for example, by selecting the minimum-norm solution.

The remainder of the article is organized as follows.

\begin{itemize}
\item Section~\ref{sec:ProblemFormul} formulates the problem in its first-order system form.

\item Section~\ref{sec:LSform} derives the nonlinear Lippmann--Schwinger integral equation of the inverse problem.

\item Section~\ref{sec:fd} presents the computation of the internal solution using a complex-symmetric Lanczos algorithm.

\item Section~\ref{sec:NumRes} reports the numerical results.

\item Section~\ref{sec:Concl} provides concluding remarks.

\item The Appendix describes the adaptive algorithm for ROM construction and the finite-difference embedding.
\end{itemize}

\section{Problem formulation}\label{sec:ProblemFormul}
\subsection{Forward and inverse problems}
We  will consider plane-wave electromagnetic wave propagation in lossy layered media in travel time coordinates $T$ on the interval $T\in[0,1]$ as derived in \cite{LossyOneDimensional},
\be
\begin{array}{rcl}
	r(T ) u(T , s) + \sigma(T)\dfrac{d}{dT} v(T, s) & = & - s  u(T, s)+ b(T), \\
	\dfrac{1}{\sigma(T)}\dfrac{d}{dT} u(T, s)  & = & - s  v(T, s), \\
	v(0)  =  0,& & u(1) = 0,
\end{array}
\label{eqn:transt}
\ee
where $b(T)$ is a localized source that approximates $\delta(T+0)$. In addition, $u$ and $v$ are,  respectively, the electric and magnetic field strength, $r$ is the electric loss function, $\sigma$ is the electromagnetic wave-impedance and $s\in\mathbb{C}^+$ is the Laplace frequency in the stable half-plane. For visco-acoustics, i.e. (damped) acoustic plane-wave propagation, we can obtain a similar first-order equation with $u$ and $v$ being, respectively, velocity and pressure using the Maxwell-dispersion model \cite{FAUCHER2023} as a loss model. Then $\sigma$ is the acoustic impedance and $r$ is the acoustic loss. 

Obviously, the variable $u$ can be eliminated from \eqref{eqn:transt}, leading to a second-order equation for $v(T)$  with a quadratic polynomial of $s$ in front of the lower-order term. We follow the derivation in \cite{LossyOneDimensional}, and after the Liouville transform $w=\sigma^{-1/2}u$, $\widehat w= \sigma^{1/2}v$ we obtain 
\be\label{eq:DiffEq}
 ({\cal A}+ s{\cal I}) \begin{bmatrix} w(T,s) \\ \widehat{w}(T,s) \end{bmatrix} = \begin{bmatrix} b(T)\\ 0 \end{bmatrix}, 
\ee
where 
\be\label{eq:operator}{\cal A}+ s{\cal I}=\calL + \calQ(T) + \calR(T) + s {\cal I}.\ee
The coefficient operator ${\cal{C}}=\calQ(T) + \calR(T)$ is divided into its symmetric components $\cal R$ corresponding to attenuation and its skew-symmetric component $\cal Q$ associated with a scattering potential. Specifically,
\be
\calQ(T) = \begin{bmatrix} 0 & \kappa(T) \\ -\kappa(T) & 0 \end{bmatrix}, 
~ \kappa(T)= \frac{d}{dT} \ln \sigma^{-\frac{1}{2}}(T) , 
~ \calR(T)= \begin{bmatrix} r(T) &0 \\ 0 & 0 \end{bmatrix}.
\ee
 Here $\calL$ is the skew-symmetric differential operator
\be
\cal{L} =  \begin{bmatrix} 0 &   \frac{d}{dT} \\  \frac{d}{dT} & 0 \end{bmatrix}.
\ee
Next, we introduce the vector-valued shorthand
\be
\phi(T) =  \begin{bmatrix} w(T) \\ \widehat{w}(T) \end{bmatrix}, \text{ and } \beta(T) =  \begin{bmatrix} b(T)\\ 0 \end{bmatrix},
\ee
such that we obtain  the PDE
\be\label{eq:PDE}
(\calL + {\cal C}(r(T),\kappa(T)) + s {\cal I} ) \phi(T,s) = \beta(T),
\ee
with boundary conditions
\be\label{eq:BC}
{
e_2^T \phi(0,s) = e_1^T \phi(1,s)=0,}
\ee
for the vector-valued field $\phi(T,s) : [0,1]\times \mathbb{C}^+ \mapsto \mathbb{C}^2$, and the complex Laplace frequency $s$ restricted to the stable half-plane $\mathbb{C}^+$. The field is excited by the source $\beta(T)$.
Henceforth, we employ the row vectors $e_1^T = [1, 0]$ and $e_2^T = [0, 1]$.

We note that all operators commute inside the $\calW$ weighted Lagrangian bilinear form 
\be
\la \phi, \psi \ra_\calW = \int_0^{1}{\rm d}T \,\, \phi(T)^T \calW\psi(T),  
\quad
\text{with} 
\quad
\calW = \begin{bmatrix}  1 &0 \\ 0 & -1\end{bmatrix}.
\ee

{\remark We also note that system of equations (\ref{eq:DiffEq}) has a port-Hamiltonian stricture with a skew-symmetric structure operator $\calL + \calQ(T)$ and a non-negative symmetric dissipation operator $\calR(T)$.  The preservation of a port-Hamiltonian structure in ROMs is a challenging problem,  and one of main difficulties we overcome in this paper is constructing a structure-preserving inversion algorithm.} 

We consider the inverse problem of recovering $r(T)$ and $\kappa(T)$ from the measurements 
\be
D(s) = \la \beta(T),  \phi(T,s) \ra_\calW, 
\label{eqn:ds}
\ee
that is, we are interested in inverting the nonlinear forward map
\be
{\cal F}: (r(T),\kappa(T)) \to D(s).
\ee

The data $D(s)$ depend non-linearly on the coefficients $r,\kappa$ since $\phi$ is unknown, and hence the inverse problem is nonlinear. 

\subsection{The rational  data model}
Since the operator ${\cal A}$ is diagonalizable (see, e.g. \cite{LossyOneDimensional}), we can write its eigendecomposition as
\be
{\cal A} = -\sum_{j=1}^\infty  \left| q_j \ra_\calW \lambda_j  \la q_j \right| + 
\sum_{j=1}^\infty\left| \overline{q_j} \ra_\calW \overline{\lambda_j}  \la \overline{q_j} \right|
\ee
with eigenfunctions $q_j(T)$ that satisfy 
\be
\la q_j(T) , q_i(T) \ra_\calW = \delta_{ij}, \text{ and }
\la q_j(T) , \overline{q_i(T)} \ra_\calW = 0,
\ee
and eigenvalues $-\lambda_j$. Here, we use the Dirac (quantum-mechanical) bracket notation in the $\calW$-outer products for the representation of operators (see, e.g. \cite{DiracNotation}).

Since we consider the inverse problem of determining $\kappa(T)$ and $r(T)$ from the measurements of a single transfer function with $\beta(T) = [b(T) , \, 0 ]^T$ as input and output, the data \eqref{eqn:ds} admit the representation.
\begin{align}\label{eqn:ds1}
D(s) &= \left<\beta(T),\phi(T)\right>_\calW = \sum_{j=1}^\infty \left(\frac{y_j}{s + \lambda_j} + \frac{\overline{y_j}}{s + \overline{\lambda_j}}\right),
\end{align}

 with the residues
 \[
 y_j = \left<\beta(T),q_j(T) \right>_\calW^2 = 
 \left( \int_{0}^1 \beta(T)^T\calW q_j(T)d T \right)^2.
 \]
 Stability and passivity of the problem yield   \[\Re( \lambda_j)\ge 0, \ \Re(y_j) \geq 0.\]
 We will obtain $[2n-1/2n]$ rational approximations (ROMs of order $2n$) \[D^{\rm ROM}(s) \approx D(s)\]  by matching either the $n$ terms closest to the origin of \eqref{eqn:ds1} or via its adaptive rational approximations, as outlined in ~\ref{ratdata}.

\section{First-Order Lippmann-Schwinger formulation}\label{sec:LSform}
Let $\phi_0(T,s)$ be a reference solution satisfying \eqref{eq:PDE} for the known $r_0(T)$ and $\kappa_0(T)$ (for example, $r_0 (T)\equiv0$ and $\kappa_0(T) \equiv 0$) with reference data $D_0(s)= \la \beta(T), \phi_0(T,s) \ra_\calW$.  Then  $\phi_0(T,s)$ satisfies equation 
\be\label{eq:DiffBackgr}
({\cal A}_0 + s{\cal I}) \phi_0(T,s) = [\calL + \calQ_0 + \calR_0 +s {\cal I}] \phi_0(T,s) = \beta(T)
\ee
with boundary conditions \eqref{eq:BC}. 

The scattered field $\phi(T,s)-\phi_0(T,s)$ now satisfies the wave equation of the background operators with forcing term depending on the field inside the unknown medium. Subtracting \eqref{eq:DiffBackgr} from \eqref{eq:DiffEq} gives
\be
({\cal A}_0 + s{\cal I}) (\phi(T,s)-\phi_0(T,s))= - \begin{bmatrix} \Delta r & \Delta \kappa \\ - \Delta \kappa & 0 \end{bmatrix} \phi(T,s)
\ee
with $\Delta r=r-r_0$ and  $\Delta \kappa= \kappa - \kappa_0$.

Then we arrive at the Lippmann-Schwinger formulation
\begin{align}\label{eq:ls}
D(s)-D_0(s) &= -\int {\rm d} T (\phi_0(T,s))^T \calW \begin{bmatrix} \Delta r & \Delta \kappa \\ - \Delta \kappa & 0 \end{bmatrix} \phi(T,s)\\
{}& =-\left<  \phi_0(T,s),  \begin{bmatrix} \Delta r & \Delta \kappa \\ - \Delta \kappa & 0 \end{bmatrix} \phi(T,s) \right>_{\calW}\nonumber,
\end{align}
which can be derived as 
\be
\begin{split}
D(s)-D_0(s)&= \la\beta(T),({\cal A}+s{\cal I})^{-1}\beta(T)\ra_\calW-\la\beta(T),({\cal A}_0+s{\cal I})^{-1}\beta(T)\ra_\calW \\
&= \la({\cal A}_0+s{\cal I})^{-1}\beta(T),({\cal A}_0-{\cal A})({\cal A}+s{\cal I})^{-1}\beta(T)\ra_\calW. 
\end{split}
\label{eq:derivation} 
\ee
 The Lippmann-Schwinger formulation allows one to preserve the block structure of the operator by explicitly imposing a zero entry in the lower right block  of the matrix on the right hand side
of \eqref{eq:ls}. This is a Fredholm integral equation of the first kind. {It is non-linear  with respect to the internal solution  $\phi$ and  the  medium  when both are unknown. However, if  one of them (the internal solution or the  medium) is known, it becomes linear with respect to the remaining parameter. 
Thus, the Born approximation (see, e.g. \cite{TarekBornRytov})
\[\phi(T,s)\approx\phi_0(T,s)\]
linearizes the Lippmann-Schwinger equation; however, it only holds for small and weak inclusions.
This Born approximation can also be iterated,  leading to an iterative Born method. Unfortunately, the Born approximation does not respect the nonlinear nature of the inverse problem and fails to recover meaningful coefficients in most scenarios.  

In this paper,  we approximate the internal field based on measured data $D(s)$. More specifically, we approximate $\phi(T,s)\approx \phi_{\rm LSL}(T,s)$ such that $\phi_{\rm LSL}(T,s)$ approximately satisfies a PDE and exactly satisfies the ROM data equation $D(s) =  \la \beta(T),  \phi_{\rm LSL}(T,s) \ra_\calW $.
Following \cite{DMZ21} we call the resulting  algorithm Lippmann--Schwinger-Lanczos (LSL) to emphasize the crucial part of the proposed approximation of the internal solution that is based on the Lanczos recursion.
By solving the LSL system using the least squares minimal norm solution, we constrain $\Delta r$ and $\Delta \kappa$, so they become regular perturbations supporting port-Hamiltonian structure by imposing positivity of $\kappa$ and non-negativity of $r$, preserving the high-frequency spectral asymptotics \cite{LossyOneDimensional}  of the background solution, thus circumventing the above-mentioned non-uniqueness of the simultaneous determination of the impedance and loss profiles. }

Through a non-linear but computationally tractable algebraic procedure, we first find an approximate transform from the eigenfunctions or (their  projections)  of the operator $\calL + {\cal C}(r(T),\kappa(T)) $ to their counterparts for  the operator $\calL + {\cal C}_0$. This allows us to expand $\phi_{\rm LSL}$ in the space of eigenfunctions (or their approximations)  of $\calL + {\cal C}_0$ while approximating $\phi$. This approximation goes beyond the Born approximation, since it takes into account the measured data. It approximately linearizes the Lippmann-Schwinger equation.


\section{Computation of the internal solution via Lanczos algorithm}\label{sec:fd}
We consider the data   model  in the form
\be
D^{\rm ROM}(s) = \sum_{j=1}^n \left( \frac{y_j}{s+\lambda_j}+\frac{\overline{y_j}}{s+\overline{\lambda_j}}\right).
\label{eqn:drom}
\ee
An obvious approximation can be obtained by truncating series  \eqref{eqn:ds1} assuming that conjugate pairs of $\lambda_j$ are ordered according to their absolute values, known as the "truncated measure"  (TM) data  \cite{LossyOneDimensional}. 
In this case, by construction $D^{\rm ROM}(s)$ is stable and passive. 
A more practically feasible alternative is adaptive rational interpolation, as presented in \ref{ratdata}. In this case $\lambda_j$ and $y_j$ are actually approximations of their counterparts in \eqref{eqn:ds1}; for simplicity, we slightly abuse notation.
\subsection{The data-driven complex symmetric Lanczos Algorithm}


The data model \eqref{eqn:drom} can be further transformed by 
applying the complex symmetric Lanczos tridiagonalization process (see Algorithm~\ref{alg:lanczos}) to the diagonal matrix of eigenvalues 
\be
\boldsymbol \Lambda=\text{diag} \left( \lambda_1,\dots,\lambda_n, \overline{\lambda_1},\dots ,,\overline{\lambda_n} \right) \in \mathbb{C}^{2n \times 2n},
\label{eqn:blambda}
\ee
and the starting vector of residues $\by = \left[ \sqrt{y_1},\dots,\sqrt{y_n},\sqrt{\overline{y_1}},\dots ,\sqrt{\overline{y_n}} \right]^T \in \mathbb{C}^{2n}$. 
Algorithm~\ref{alg:lanczos} computes a tridiagonal complex-symmetric matrix 
$\bT \in \mathbb{C}^{2n \times 2n}$ and a matrix
$\bV = [\bv_1, \ldots, \bv_{2n}] \in \mathbb{C}^{2n \times 2n}$ satisfying Lanczos identities
\be
\bT=\bV^T \boldsymbol \Lambda \bV, 
\quad \bV^T\bV=\bI_{2n}, \text{~and~}
\bv_1 = \frac{1}{\sqrt{\by^{T} \by}} \by,
\label{eqn:lanczosid}
\ee
that allow us to rewrite the data model as a quadratic form
\be
\label{eq:TF_tridiag_2}
D^{\rm ROM}(s)= \by^T (\boldsymbol \Lambda + s \bI_{2n})^{-1} \by=\frac{1}{\widehat\gamma_{1}} \bee_{1}^{T} (\bT+s{\bI_{2n}})^{-1} \bee_{1},
\ee
where $\bee_{1} \in \mathbb{C}^{2n}$ is the first column of the $2n \times 2n$ identity matrix $\bI_{2n}$.

We observe that the nonzero real part of the eigenvalues in \eqref{eqn:blambda} gives rise to a nonzero main diagonal in the tridiagonal matrix $\bT$. Specifically, the off-diagonal entries $\beta_{j+1}=\bV_{j+1}^T \boldsymbol \Lambda \bV_j$ are purely imaginary and the diagonal entries $\alpha_j=\bV_j^T \boldsymbol \Lambda \bV_j$ are purely real due to the conjugate structure of $\bV$ and $\boldsymbol \Lambda$. 

\begin{algorithm}
\caption{The Lanczos process for complex symmetric matrices.}
\label{alg:lanczos}
\begin{spacing}{0.9}
\begin{algorithmic}[1]
\State Normalize: $\bv_1 = \by/\sqrt{\by^{T} \by}$
\State $\bw = \boldsymbol\Lambda\bv_1$ 
\State $\alpha_1 = \bw^T \bv_1$
\State $\bw = \bw - \alpha_1 \bv_1$
\For{$j = 2, \ldots, 2n$ }
\State $ \beta_j = \sqrt{\bw^T \bw}$ 
\If{$\beta_j \neq 0$}
\State $\bv_j = \bw/\beta_j$
\State $\bw = \boldsymbol\Lambda \bv_j$ 
\State $\alpha_j = \bw^T \bv_j$ 
\State $ \bw = \bw - \alpha_j \bv_j - \beta_j \bv_{j-1}$
\Else \,\,Breakdown
\EndIf
\EndFor
\end{algorithmic}
\end{spacing}
\end{algorithm}

{\remark An ambiguity arises in choosing $\beta_j$ at step 6. To make it consistent for background and true medium, we always choose the principal branch of square root and consider $\bw^T \bw+\epsilon$ with $\epsilon>0$ if $\bw^T \bw$ is on the branch cut, i.e., negative real.}

The inverse problem can be solved using the finite-difference embedding of $\bT$  \cite{LossyOneDimensional}. However, this approach requires constrained nonlinear optimization even for moderate variation of losses and, unlike the LSL, does not have obvious extension to  multidimensional problems, even for lossless media. We outline the finite-difference embedding for completeness  in \ref{app:B}.

\subsection{Computation of internal solution}

To compute the approximate internal solution, we collect the first $n$ eigenfunctions and their complex conjugates in a quasi-matrix
\be
Q = [q_1, q_2, \dots, q_n, \overline{q}_1,\overline{q}_2, \ldots, \overline{q}_n], 
\ee
with the corresponding eigenvalues on the diagonal of $\boldsymbol\Lambda$ defined in \eqref{eqn:blambda}. 
We obtain
\begin{align}\label{eq:eigdecomp}
\phi(T,s)&= ({\cal A}+s{\cal I})^{-1}\beta(T)\\
{}& = {Q} (\boldsymbol \Lambda + s\mathbf{I}_{2n})^{-1} \la {Q} | \beta(T) \ra_\calW
.\end{align}
where $\la {Q} | \beta(T) \ra_\calW=\begin{bmatrix} \la q_1,\beta(T)\ra_\calW \\ \vdots \\ \la \bar{q}_n,\beta(T)\ra_\calW \end{bmatrix}$.

We approximate $\phi(T,s)$ using Lanczos-ROM state solutions. 
Applying Algorithm~\ref{alg:lanczos} to the poles and residues of $D^{\rm ROM}$, and assuming that it completes all $2n$ steps without breakdown, we obtain the matrices $\bT,\bV\in\mathbb{ C}^{2n\times 2n}$. Provided that the forcing function $\beta(T)$ satisfies \eqref{eqTMb}, Algorithm~\ref{alg:lanczos} performs an exact tridiagonalization. That is, in the coordinates of the eigenfunctions, using \eqref{eqhath1}, we obtain the Lanczos identity $\boldsymbol\Lambda\bV=\bV\bT$, where the first vector of the Lanczos basis can be expressed in terms of residues as
\be
  \bv_1= \frac{1}{\sqrt{\sum_{j=1}^n(y_j+{\overline y_j})}}[\sqrt{y_1},\ldots, \sqrt{y_n},\sqrt{\overline y_1},\ldots, \sqrt{\overline y_n}]^T.
 \ee
 Alternatively, the $\cal W$-bilinear form of the forcing function $\beta(T)$ allows us to write
\be
   \bv_1 =({\la\beta(T),\beta(T)\ra_{\cal W}})^{-\frac12}\la {Q} | \beta(T) \ra_\calW
   \label{eqn:bv1beta}
\ee
Denoting $\la\beta(T),\beta(T)\ra_{\cal W} =\int_0^1 b(T)^2 dT = \|b(T)\|^2$, the Lanczos identity can be written in the state space coordinates as 
  \be \label{eq:LanczRec} {\cal A}{Q}\bV={Q}\bV\bT, \qquad 
  {Q}\bv_1= \frac{1}{\|b(T)\|}\beta(T).\ee
Alternatively, the internal field $\phi(T, s)$ from \eqref{eq:eigdecomp} can be expressed as
\be\label{eq:Lancz}
 \phi(T,s)  = {Q} \bV (\bT + s\mathbf{I}_{2n})^{-1} \bV^T \la {Q} | \beta(T) \ra_\calW={Q} \bV (\bT + s\mathbf{I}_{2n})^{-1}\|b(T)\| \mathbf{e}_1,
 \ee
where the identity $ \bV^T \la Q | \beta(T) \ra_\calW=\|b(T)\| \mathbf{e}_1$ follows from \eqref{eqn:bv1beta} and \eqref{eqn:lanczosid}. 

At the heart of our approach lies the crucial observation (as discussed in the next section) that the Lanczos basis in travel time coordinates ${Q} \bV$ is only weakly dependent on the unknown coefficients $\kappa(T)$ and  $r(T)$. Thus, it can be replaced with ${Q_0} \bV_0$, computed for a background with known coefficients $\kappa_0(T)$ and $r_0(T)$. That is, in the LSL algorithm, we approximate the internal solution via \be\label{eq:lancinv}
{Q} \bV \approx {Q}_0 \bV_0,\ee 
and 
\be\label{eq:intsol}
\phi(T,s)  \approx\phi_{\rm LSL}\equiv {Q}_0 \bV_0 (\mathbf{T} + s\mathbf{I}_{2n})^{-1} \|b(T)\| \mathbf{e}_1.
\ee
Substitution of the data-driven estimate \eqref{eq:intsol} into the Lippmann-Schwinger equation \eqref{eq:ls}  effectively linearizes it with respect to unknown coefficients $\kappa$ and $r$.

\begin{remark}
We note that for the TM case  $D^{\rm ROM}(s)\equiv D^{\rm ROM}(s)$ if and only if  the forcing function $\beta(T)$ admits a finite expansion in the eigenfunction basis:
\be
\label{eqTMb} 
\beta(T)=\sum_{j=1}^n \left( e_1^T q_j(0) q_j + 
\overline{e_1^T q_j(0) q_j} \right).
\ee
In the inverse problem, we use $\beta(T)=[\delta(T), 0]^T$ or its smooth approximation. The eigenfunctions $q_j(T)$ depend on the unknown PDE coefficients, and so does $\beta(T)$ given by \eqref{eqTMb}. However, for large $j$  the eigenfunctions become asymptotically close to the case $\calQ(T)=0,  \calR(T)=0$ \cite{LossyOneDimensional}. Thus, we assume that $n$ is large enough so that $\beta(T)$ is well approximated in the eigenfunction bases of both the unknown medium and the background medium.
\end{remark}

 \subsection{Qualitative justification of the internal solution estimate}
  \label{sect:explanation}

	The data-driven ROM inversion based on  approximation \eqref{eq:lancinv}. This effect was explained for the lossless (Hamiltonian) case with $r(t)\equiv 0$ in \cite{DMTZ16,BDMMZ2020,BorceaZimmerling}, applied in \cite{baker25, Abilgazy}, and analyzed in \cite{NormG,druskin2025optimality}. 

From \eqref{eq:lancinv} we define the transmutation operator $\mathbf{R}=\mathbf{V_0}\mathbf{V}^{-1}$ mapping background eigenfunctions to approximate true ones
\be\label{eq:trans} \mathbf{Q} \approx \mathbf{Q}_0 \mathbf{R}. \ee

For discrete, lossless, first-order problems, such an operator  was introduced in \cite{borcea2005continuum}, following an unpublished note by Natterer.  In conventional second-order formulations, continuous transmutations operators play a central role in the Marchenko–Gelfand–Levitan inverse scattering theory (see, e.g. \cite{Kravchenko2020SomeRD}). In \cite{druskin2025optimality}, it was shown that $\mathbf{R}$ yields convergent approximations to internal solutions of the second-order problem in the time domain as the ROM dimension tends to infinity.  We therefore conjecture similar convergence of \eqref{eq:lancinv}.

However, we should note that the first-order formulation differs essentially from the second-order one. In particular, \eqref{eq:trans} involves  conjugate pairs of eigenfunctions, corresponding to both downward- and upward-propagating waves. Therefore, to extend the analysis of  \cite{druskin2025optimality} to the first-order case, one needs to include  both  types of waves in the transmutation operator as was done in the original work \cite{borcea2005continuum} using the lossless first-order formulation. Otherwise, convergence to the true first-order solution can be lost, as in  \cite{borcea2025reduced}. 

Even with both wave-types included, extending the proof of \cite{druskin2025optimality} is obstructed by possible breakdowns of the non-Hermitian Lanczos method. Therefore, here we limit ourselves to qualitative justification.

	Let us first consider the formulation in the time-domain (see \cite{DMTZ16,BorceaZimmerling}). It is based on data-driven orthogonalization of a causal solution, transforming them into the shift-invariant waves corresponding to the homogeneous medium. However, as mentioned in the introduction, the time-domain formulation with finite time observation does not resolve between lossy media and lossless perturbations with trapped modes. In the bounded domains losses manifest themselves by real positive components of spectra, which can be observed only in the frequency domain, or equivalently,   requires sampling the semi-infinite time interval.  Algebraically, the failure of the time-domain approach in the dissipative case is  due to the well-known fact that the kinetic energy of the shift-invariant solutions is equal to the potential one. This annuls the  Lagrangian bi-linear form and  leads to the breakdown (or near breakdown) of the   non-Hermitian Lanczos algorithm, known in the computational linear algebra community \cite{Freund93}.
	
	Thus, the remaining option for our case is the alternative approach, which exploits the Bubnov-Galerkin equivalent of the finite-difference ROM realization (see, e.g. \cite{BDMMZ2020}). The simplest approach to compute the data-driven solution is to use the discretization on the TM grid. As discussed in Section~\ref{sec:fd}, the finite-difference discretization \eqref{eqn:romSISO} yields convergence of the discrete coefficients to the true continuum ones in the sense of Lagrangian bilinear form. The staggered finite-difference scheme on the TM grid has a first-order approximation in the interior and an exact transfer function  \cite{borcea2005continuum}. This is rather exceptional, as finite-difference Gaussian rules generally provide second-order approximation at the interior and spectral convergence of the transfer function. Therefore, using \eqref{eqn:romSISO} with the inverted coefficients, as described in Section~\ref{sec:fd}, gives the first-order approximation of  the solution of \eqref{eqn:transt}. After discrete Liouville transform
	$w_i=\frac{1}{\sqrt{\sigma_i}} u_i$, $\widehat w_i= \sqrt{\widehat{\sigma}_i} v_i$ we obtain the first-order approximation to the solution of \eqref{eq:DiffEq}, or equivalently, at the grid nodes
	\be\label{eq:FD}
	  (\bT + s\mathbf{I}_{2n})^{-1}\|b(T)\|\mathbf{e}_1.
	\ee
	The approximation \eqref{eq:FD} is purely data-driven (it  only requires $\bT$ which is computed from the data) and in principle its interpolation can be used in 
	\eqref{eq:Lancz}. However, the first-order approximation puts it at a disadvantage compared to the exact solution given by \eqref{eq:ls} for the case of small perturbations, when the Born approximation is accurate enough. To circumvent this disadvantage, let us compare \eqref{eq:ls} and \eqref{eq:FD} in the grid nodes. They approximate the same function $\phi(T,s)$ but differ by operator-valued factor $\mathbf{Q} \bV$, so this matrix must be close to identity at the grid nodes.
	 The same reasoning holds for $\mathbf{Q}_0\bV_0$, so one can replace the other, thus justifying \eqref{eq:intsol}. Since the eigenfunctions $q_i$ become asymptotically independent on $\kappa$ and $r$ as $n\rightarrow \infty$ \cite{BDK2005}, the range of 
	 $\mathbf{Q} \bV$ becomes just close to one of the eigenspaces $\mathbf{Q}_0$, and it can be used for interpolation between the grid nodes. Indeed, this qualitative reasoning cannot replace rigorous analysis that will be the subject of a forthcoming work.

\section{Numerical Results}\label{sec:NumRes}
\subsection{Spectral measure ROM}
To show that the LSL approach is superior to the embedding strategies used in \cite{LossyOneDimensional}, 
we rerun the same numerical example where the direct embedding starts to break down. 
The potential $\kappa$ used in this manuscript is related to the impedance $\sigma$ from \cite{LossyOneDimensional} 
via $\kappa(T) = \frac{d}{dT} \ln \sigma^{-\frac{1}{2}}$.  

\begin{figure}[ht]
\centering
\renewcommand{\arraystretch}{1.2} 
\begin{tabular}{m{0.5cm} c c} 
   & \textbf{Imaginary part} & \textbf{Real part} \\[0.5em]
\rotatebox[origin=c, y=7cm]{90}{\textbf{Primary wave $w$}}
   & \includegraphics[height=4.67cm]{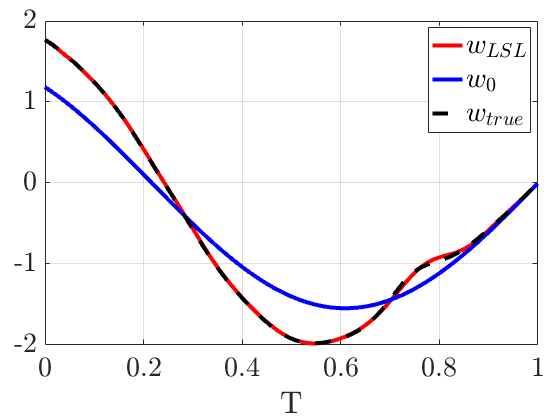} 
   & \includegraphics[height=4.65cm]{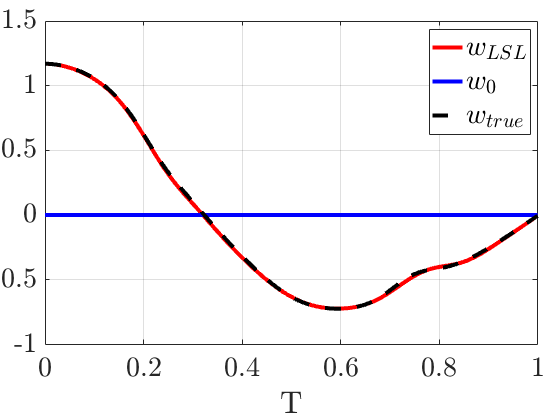} \\[-4cm]
\rotatebox[origin=c, y=7cm]{90}{\textbf{Dual wave $\widehat w$}}
   & \includegraphics[height=4.65cm]{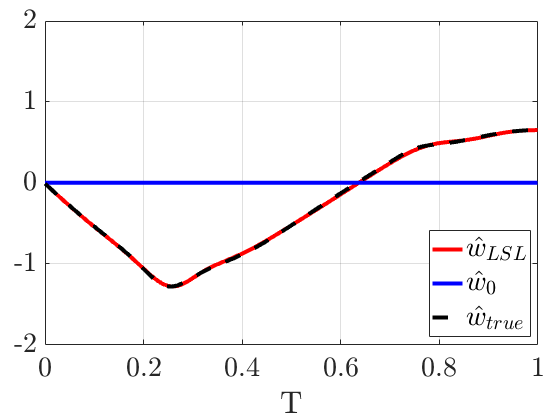} 
   & \includegraphics[height=4.65cm]{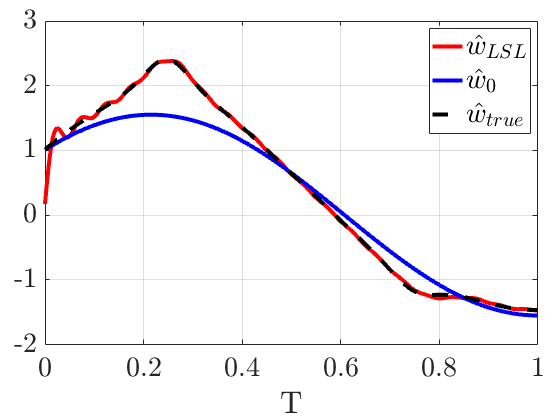} \\[-4cm]
\end{tabular}
\caption{Internal solutions for $n=40$ at frequency $s = 4i$. 
The top row shows the primary wave $w$, and the bottom row shows the dual wave $\widehat{w}$. 
The left column displays the imaginary part, while the right column displays the real part. 
Note that the chosen forcing function is equivalent to enforcing the condition $\lim_{T\to 0^+} \widehat{w}(T) = 1$. 
The wave of the background $w_0$, $\widehat{w}_0$ is the Born approximation frequently used in the Lippmann--Schwinger equations.}
\label{fig:internal_solutions}
\end{figure}

Figure~\ref{fig:internal_solutions} illustrates the internal solutions for $n=40$ at the frequency $s = 4i$, 
showing both the primary wave $w$ and the dual wave $\widehat{w}$. 
The background solution has a purely imaginary primary wave and a purely real dual wave. 
Clearly, the internal solution found via the LSL is more accurate than the Born solution.  

\begin{figure}[ht]
\centering
\includegraphics[height=5cm]{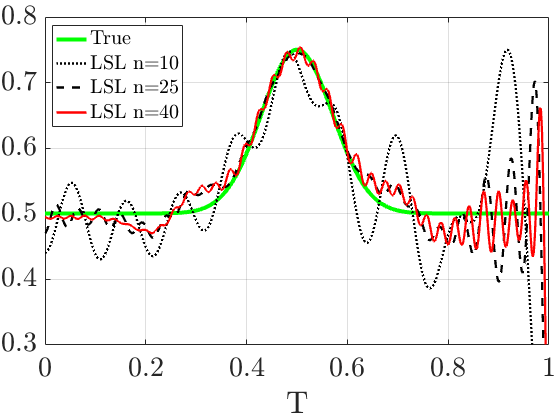}
\hfil
\includegraphics[height=5cm]{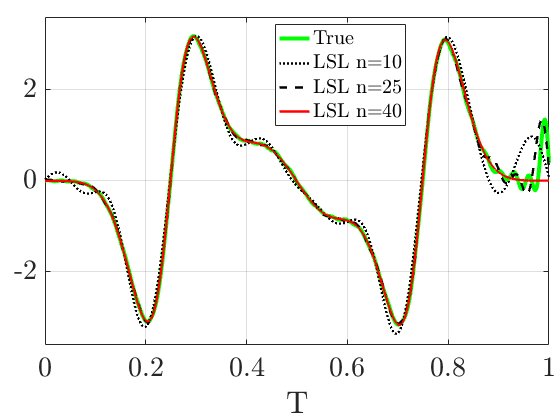}
\caption{Inversion results after a \textbf{single} iteration of the LSL algorithm for $n=10,25,$ and $40$, with Tikhonov regularization. Left: recovery of a large attenuation profile $r(T)$. Right: First order potential $\kappa(T)$.}
\label{fig:inversion_results}
\end{figure}

These solutions form the basis for the inversion shown in Figure~\ref{fig:inversion_results}, 
which presents the reconstructions obtained after a single iteration of the LSL algorithm for $n = 10, 25$, and $40$. 
The respective eigenvalues with the largest imaginary part for these three inversions are  
$\lambda^{max}_{n=10} = 0.2719 + 29.8041i$,  
$\lambda^{max}_{n=25} = 0.2721 + 76.6789i$, and  
$\lambda^{max}_{n=40} = 0.272 + 123.41i$.  
The potential is well recovered for all settings of $n$. 
The loss function is recovered well for larger values of $n$, although the recovery exhibits some oscillation and would benefit from gradient penalization. The high losses make it difficult to accurately recover the medium at the end of the domain; however, in the rest of the domain, we have quantitative agreement after a single iteration of the LSL method.

\subsection{Data-driven adaptive ROM}
\begin{figure}[H]
\centering
\renewcommand{\arraystretch}{1.2} 
\begin{tabular}{m{0.5cm} c c} 
   & \textbf{Losses} & \textbf{Impedance} \\[0.5em]
\rotatebox[origin=c, y=7cm]{90}{\textbf{0\% noise}}
   & \includegraphics[height=4.75cm]{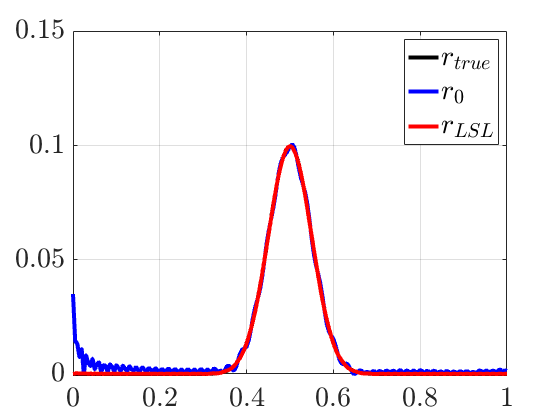} 
   & \includegraphics[height=4.75cm]{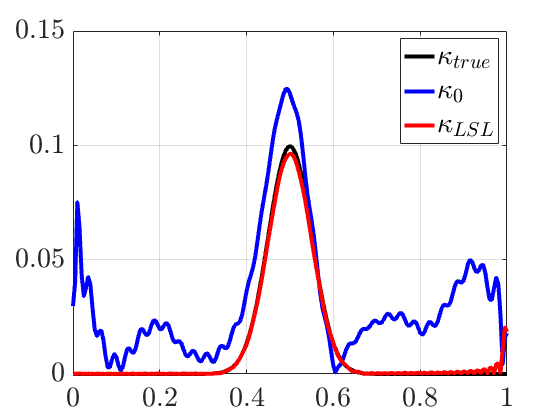} \\[-4cm]
\rotatebox[origin=c, y=7cm]{90}{\textbf{20\% noise}}
   & \includegraphics[height=4.75cm]{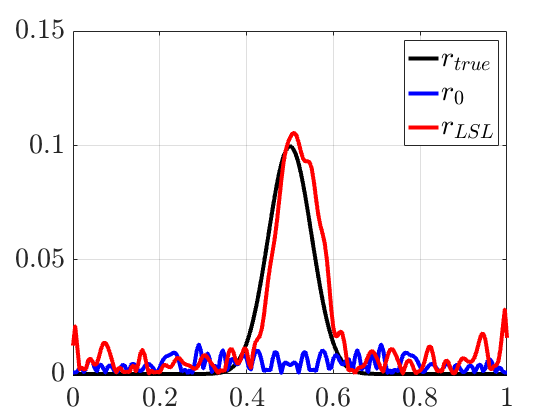} 
   & \includegraphics[height=4.75cm]{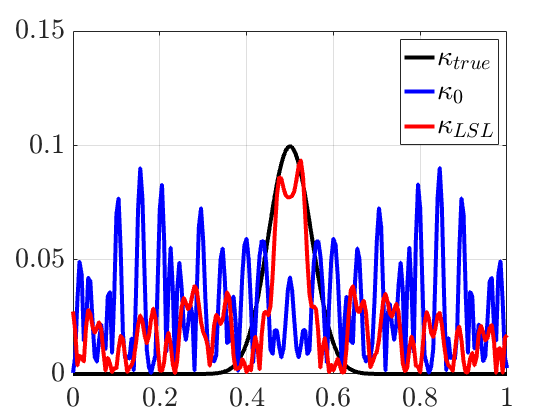} \\[-4cm]
\end{tabular}
\caption{Numerical examples with losses and impedance profiles given by Gaussians (black curves). In the noiseless case (top plots), the LSL approach (red curves) allows to produce qualitatively better reconstructions compared to the Born method (blue curves). Bottom plots show that LSL remains stable even with 20\% noise whereas the Born method fails. }
\label{fig:adapt}
\end{figure}

In the next set of numerical experiments, we compared LSL with the adaptive ROM to construct internal solutions against Born approach. We considered problem (\ref{eq:PDE}) with loss and impedance profiles shown in black in Fig.\ref{fig:adapt}. Data for 15000 frequencies $\omega_j$, uniformly distributed in the frequency interval $[0.1;100]$, were generated through numerical modeling with the second-order finite-difference scheme. The adaptive approach required 182 frequencies to achieve a $10^{-13}$ data error $\|\widehat{D}-D\|_\infty$ for the data without noise. We plugged the obtained data-driven ROM into (\ref{eq:ls}) discretized using quadrature with 1000 nodes. As can be seen, even in the noiseless scenario (top plots in Fig. \ref{fig:adapt}), the LSL method (red curves) outperforms the Born approach (blue curves) for the chosen frequencies. Moreover, adding significant $20\%$ noise (bottom plots in (Fig. \ref{fig:adapt})) still allows the LSL approach to produce reasonable results, whereas the Born approach completely breaks down. 

\vspace{0.5cm}
\captionsetup{labelfont=bf} 

\captionsetup{labelfont=bf} 

\vspace{0.5cm}
\captionsetup{labelfont=bf} 



\section{Conclusions}\label{sec:Concl}
We have presented an extension of the Lippmann–Schwinger–Lanczos (LSL) algorithm for the simultaneous reconstruction of loss and impedance profiles in one-dimensional media. The method is formulated through a first-order nonlinear Lippmann–Schwinger integral equation and employs a data-driven approximation of the internal field within the unknown region using the non-Hermitian, complex-symmetric Lanczos algorithm. Although the latter may, in theory, exhibit algorithmic breakdowns—none of which were observed in our experiments—this limits a fully rigorous extension of the justification known for the lossless case. Nevertheless, we provide intuitive reasoning for the validity of the approach via finite-difference embedding. Numerical experiments demonstrate its efficiency and robustness, even in the presence of significant measurement noise. Future work will focus on extending this framework to multidimensional problems, with particular emphasis on applications in ground-penetrating radar imaging in complex environments.
\section*{Acknowledgements}
V. Druskin was partially supported by AFOSR grants FA 9550-20-1-0079, FA9550-23-1-0220,  and NSF grant  DMS-2110773. A.V. Mamonov was supported in part by the U.S. National Science Foundation under award: DMS-2309197. This material is based on the research supported in part by the U.S. Office of Naval Research under award number N00014-21-1-2370 to A.V. Mamonov. M. Zaslavsky was partially supported with AFOSR grant  FA9550-23-1-0220.

\bibliographystyle{elsarticle-num} 
\bibliography{bib}


\appendix

\section{Adaptive data-driven ROM}\label{ratdata}
\label{sec:adaptiveROM}
Although the truncated measure data model is the simplest, it is also not the most realistic. A more realistic data model can be obtained if we consider the data-driven Bubnov-Galerkin projection framework. 

Consider Bubnov-Galerkin projection of (\ref{eq:operator}) onto the rational Krylov subspace spanned by the entries of 
\be
U=[\phi_1=\phi(i\omega_1),\ldots, \phi_{n}=\phi(i\omega_n)],
\ee
with respect to $\calW$ bilinear form. Here $\phi(i\omega_j)$, $j=1,\ldots,n$, are the solutions of \eqref{eq:operator} for $s=i\omega_j$, $j=1,\ldots,n$. 
For the remainder of this work, we assume that $n$ is even and the frequencies $\omega_j$, $j=1,\ldots ,n$ come in positive and negative pairs, i.e. $\omega_{2k}=-\omega_{2k-1}$,  $k=1,\ldots ,n/2$. Then, the projected system (\ref{eq:operator}) takes the form
\be
\label{eq:gal}
\bS\bphi+s\bM\bphi=\bb,
\ee
where $\phi\approx U\bphi$, $\bS \in \CC^{n \times n}$ is a symmetric indefinite stiffness matrix with entries 
\be 
S_{pq} = \left<\phi_p(T),{\cal A}\phi_q(T)\right>_\calW, 
\quad p,q=1,\ldots,n,
\ee 
while $\bM \in \CC^{n \times n}$ is a symmetric indefinite mass matrix with entries 
\be
M_{pq}=\left<\phi_p(T),\phi_q(T)\right>_\calW,
\quad p,q=1,\ldots,n,
\ee 
and the entries of $\bb \in \CC^n$ are given by 
\be
b_q=\left<\beta(T),\phi_q(T)\right>_\calW = 
D(i\omega_q), \quad q=1,\ldots,n.
\ee 
The transfer function of the Bubnov-Galerkin projected system is given by
\be
\label{romdata}
D^{\rm ROM}(s) = \bb^T\bphi(s),
\ee
and it is related to the true data via the matching conditions
\be
D_q = D(i \omega_q) = D^{\rm ROM}(i \omega_q), ~
D^\prime_q = \frac{d D}{ds} (i \omega_q) = \frac{d D^{\rm ROM}}{ds} (i \omega_q), ~ q = 1,\ldots,n.
\ee
We note that although internal solutions $\phi_p$, $p=1,\ldots n$ are not accessible since $\calQ$ and $\calR$ in \eqref{eq:operator} are unknown, the entries of the mass and stiffness matrices can still be computed directly from the data using the Loewner framework \cite{Antoulas}. 
Explicitly,
\begin{eqnarray}
M_{pq} & = & \dfrac{D_p-D_q}{i\omega_p-i\omega_q}, \quad p \neq q, \label{massmtr} \\
M_{pp} & = & - D^\prime_p,  \label{massmtrdiag} \\
S_{pq} & = & \dfrac{\omega_p D_q - \omega_q D_p}{\omega_p -\omega_q}, \quad p \neq q, \label{stifmtr} \\
S_{pp} & = & D_p  + i \omega_p D^\prime_p. \label{stifmtrdiag}
\end{eqnarray} 


The choice of interpolation points $\omega_1,\ldots ,\omega_n$ is crucial for the performance of the constructed ROM. For diffusive problems, the priori selection of interpolation frequencies was proposed in \cite{dkz} and is based on the solution of Zolotarev problem. The approach was then generalized for high and infinite-order self-adjoint dynamical systems \cite{dz}. In \cite{dlz,druskin2011adaptive}, an adaptive approach was developed to select interpolation points in a greedy manner for diffusive problems. This method outperforms a priori selection, particularly for cases of strongly non-uniform spectral measures. In this paper, we extend the approach to handle damped wave problems. The approach is summarized in Algorithm~\ref{alg:adpol}. We also note that there it is partly similar to AAA algorithm \cite{AAA}.

\begin{algorithm}[ht]
\caption{Adaptive rational fitting procedure for data.}
\label{alg:adpol}
\begin{algorithmic}[1]
\State For the given frequency range of interest $[\omega_{min},\omega_{max}]$, set $\omega_1=\sqrt{\omega_{min}\omega_{max}}$, and $\omega_2=-\omega_1$,  
\For{$n = 2, 4, 6, \ldots$ }
\State Compute the matrices $\bS$ and $\bM$ via \eqref{massmtr}--\eqref{stifmtrdiag} for interpolation frequencies $\omega_1,\ldots ,\omega_n$;
\State Compute ROM data $D^{\rm ROM}$ for $\omega\in[\omega_{min},\omega_{max}]$ via \eqref{eq:gal} and \eqref{romdata};
\State Compute the next interpolation frequency by maximizing the error:
\be
\omega_{n+1} = \mathop{\text{argmax}}\limits_{\omega \in [\omega_{min},\omega_{max}]} \left| D^{\rm ROM}(i\omega) - D(i \omega) \right|
\ee
If there are several frequencies for which the maximum is attained, select any one such frequency. 
\State Set $\omega_{n+2}=-\omega_{n+1}$.
\EndFor
\end{algorithmic}
The poles and residues of obtained ROM in form \eqref{eqn:drom} are computed by solving the eigenproblem for matrix pencil  $(\bS,\bM)$.
\end{algorithm}
Note that the obtained ROM is not structure-preserving, meaning it may not inherit passivity or even stability of the original system \eqref{eq:operator}.

For the cases when $\bM$ is ill-conditioned Algorithm~\ref{alg:adpol} requires a regularization of matrix pencil $(\bS;\bM)$ similar to the one described in \cite{baker25, Abilgazy}.

{\color{blue}}
\section{Finite-difference embedding }\label{app:B}

We follow  \cite{LossyOneDimensional} and consider the TM data.
We  re-parametrize the entries of $\bT$ as
\be
\label{eq:Tgeneric}
\bT=
\setlength{\arraycolsep}{4pt}
\begin{bmatrix}
r_1 & \dfrac{-i}{\sqrt{\gamma_{1} \widehat\gamma_{1}}} & 0 & &0 \\[-5 pt]
\dfrac{-i}{\sqrt{\gamma_{1} \widehat\gamma_{1}}}  & \widehat r_1 & \dfrac{i}{\sqrt{\gamma_{1} \widehat\gamma_{2}}} &0& \\[-5 pt]
\vdots & \dfrac{i}{\sqrt{\gamma_{1} \widehat\gamma_{2}}} & \ddots & \ddots\\[-5 pt]
\vdots&     & \ddots & r_n & \dfrac{-i}{\sqrt{\gamma_{n} \widehat\gamma_{n}}} \\[-5 pt]
0 & \cdots &  0 & \dfrac{-i}{\sqrt{\gamma_{n} \widehat\gamma_{n}}} & \widehat r_n\\[-5 pt]
\end{bmatrix}.
\ee
This reparameterization allows us to directly interpret the matrix $\bT$ as finite-difference discretization of a related PDE, which was used to formulate a finite-difference inversion approach in \cite{LossyOneDimensional}. We briefly review this approach and highlight the issue caused by dual losses $\widehat r_j$ in this parameterization. In the LSL approach, the same tridiagonal matrix is used to generate internal solutions, which circumvents some of the problems arising from the direct embedding of the tridiagonal matrix.
The matrix $\bT$ in \eqref{eq:Tgeneric} can be viewed as a symmetrized finite-difference operator with real coefficients mimicking the discretization of \eqref{eqn:transt}, see \cite{LossyOneDimensional}. Explicitly,
\be
\begin{array}{rcl}
	{r_j}u_j  +  \sigma_j \dfrac{v_{j} - v_{j-1}}{\widehat{ h}_j} & = &  -s u_j, \\
	\dfrac{u_{j+1} - u_{j}}{\widehat{\sigma}_j h_j} + \widehat{r}_j  v_j & = & - s  v_j, \\
	v_0 & =& -1, \quad u_{n+1} = 0.
\end{array}
\label{eqn:romSISO}
\ee
where $h_i$ and $\widehat{h}_i$ are, respectively, the primary and dual steps of a staggered
finite-difference scheme, $\sigma_j$ and $\widehat{\sigma}_j$ are the discrete impedance defined, respectively, on the primary and dual grid nodes. Similarly, the primary and dual dampers are $r_j$ and $\widehat{r}_j$. 
These $4n$ independent real parameters are extracted from knowledge of $2n$ complex parameters $\{(y_j,\lambda_j)\}_{j=1}^{n}$ and the condition 
\be\label{eqhath1} \|\mathbf{\beta}\|=\sqrt{\sum_{j=1}^n(y_j+{\overline y_j})}.\ee Further, $\bT$ can be transformed to the same number of parameters
 $r_j$,  $\widehat r_j$, $\dfrac{\sigma_j}{\widehat{ h}_j}$ and ${\widehat\sigma_j}{ h}_j$, $j=1,\ldots, n$ via a simple one-to-one recursive algorithm.
For the lossless ($r\equiv 0$) case, if $\bT$ is computed from truncated measure data $D^{\rm ROM}$ then $h_i$ and $\widehat{h}_i$  are, respectively, the primary and dual steps of a so-called truncated measure (TM) grid, which gives (under some regularity conditions) pointwise convergence of $\sigma_i$ and $\widehat{\sigma}_i$ to the true medium parameters on the grid nodes with steps (dual) $\widehat h_i$ and (primary) $h_i$, respectively.  The TM grid can be computed for any known regularly enough impedance distribution, e.g. $\sigma\equiv 1$. It was shown in \cite{borcea2005continuum} that such an approach gives convergence to the continuum case as $n\rightarrow \infty$, which is equivalent to the fact that the TM grids for all regular enough $\sigma$ become asymptotically close to each other for large $n$. The TM grid is a particular case of the so-called finite-difference Gaussian rules or spectrally matched grids \cite{DK99}.

The same approach can be extended straightforwardly for the lossy case with $r(x)=const$, in which case $r_i$
approximates that constant, with the grid computed for the lossless case with known $\sigma(x)$. The situation becomes more complicated for the lossy case with variable $r$. The Lanczos algorithm cannot be constrained, which recursively maps $4n$ parameters of $\eqref{eqn:drom}$ to
the same number of parameters in $\|\beta(T)\|$ and $\mathbf{T}$. That results in the appearance of non-physical dual losses $\widehat r_i$ that can even become negative, even though the measurements have a passive response. However, as shown in \cite{LossyOneDimensional}, for slowly varying $r(T)$, it can be approximated by the interpolated difference of $r_i$ and $\widehat r_i$, which is related to the appearance of $\pm 1$ in the weight matrix $\calW$ of the Lagrangian bilinear form. Unfortunately, this approximation loses accuracy for a strong variation of $r(x)$ that requires a constrained (nonconvex)  misfit optimization with multiple solutions of the forward problem.

\end{document}